\documentclass[onefignum,onetabnum]{siamonline220329}

\usepackage{epstopdf}
\usepackage{algorithmic}
\ifpdf
  \DeclareGraphicsExtensions{.eps,.pdf,.png,.jpg}
\else
  \DeclareGraphicsExtensions{.eps}
\fi

\newsiamthm{assumption}{Assumption}

\headers{Bayesian Inference with Projected Densities}{J.M. Everink, Y. Dong, M.S. Andersen}

\title{Bayesian Inference with Projected Densities\thanks{\funding{This work was supported by The Villum Foundation (grant no. 25893).}}}

\author{Jasper M. Everink\thanks{Department of Applied Mathematics and Computer Science, Technical University of Denmark. Richard Petersens
Plads, Building 324, DK-2800 Kgs. Lyngby, Denmark. (\href{mailto:jmev@dtu.dk}{jmev@dtu.dk}, \href{mailto:yido@dtu.dk}{yido@dtu.dk}, \href{mailto:mskan@dtu.dk}{mskan@dtu.dk})}\and Yiqiu Dong\footnotemark[2] \and Martin S. Andersen\footnotemark[2]}

\usepackage{float}
\usepackage{bm}

\usepackage{amsmath,amssymb,amsopn,amsfonts, mathabx}

\usepackage{caption}
\usepackage{subcaption}
\usepackage{newfloat}

\DeclareFloatingEnvironment[
    fileext=los,
    listname={Algorithm},
    name=Sampler,
    within=section,
]{algo}

\usepackage{amsopn}

\DeclareMathOperator{\bd}{bd}
\DeclareMathOperator{\relint}{relint}
\DeclareMathOperator{\relbd}{relbd}
\DeclareMathOperator{\aff}{aff}
\DeclareMathOperator*{\argmin}{argmin}

\newcommand{\erfc}{\mathop{\mathrm{erfc}}}

\newcommand{\nullspace}{{\text{Null}}}
\newcommand{\reals}{\mathbb{R}}
\newcommand{\symm}{\mathbb{S}}

\renewcommand{\vec}[1]{\bm{#1}}
\newcommand{\set}[1]{\mathbf{#1}}

\usepackage{tikz}
\usetikzlibrary{patterns}
\usepackage{graphicx}
\usepackage{hyperref}

\makeatletter
\newcommand*{\addFileDependency}[1]{
  \typeout{(#1)}
  \@addtofilelist{#1}
  \IfFileExists{#1}{}{\typeout{No file #1.}}
}
\makeatother

\ifpdf
\hypersetup{
  pdftitle={Bayesian Inference with Projected Densities},
  pdfauthor={J.M. Everink, Y. Dong, M. S. Andersen}
}
\fi

\begin{document}

\maketitle

\begin{abstract}
Constraints are a natural choice for prior information in Bayesian inference. In various applications, the parameters of interest lie on the boundary of the constraint set. In this paper, we use a method that implicitly defines a constrained prior such that the posterior assigns positive probability to the boundary of the constraint set. We show that by projecting posterior mass onto the constraint set, we obtain a new posterior with a rich probabilistic structure on the boundary of that set. If the original posterior is a Gaussian, then such a projection can be done efficiently. We apply the method to Bayesian linear inverse problems, in which case samples can be obtained by repeatedly solving constrained least squares problems, similar to a MAP estimate, but with perturbations in the data. When combined into a Bayesian hierarchical model and the constraint set is a polyhedral cone, we can derive a Gibbs sampler to efficiently sample from the hierarchical model. To show the effect of projecting the posterior, we applied the method to deblurring and computed tomography examples.
\end{abstract}



\section{Introduction}

The goal of Bayesian inference is to analyze the probability distribution of a parameter obtained through Bayes' theorem \cite{gelman1995bayesian}. Bayes' theorem states that an initial distribution $\pi(\vec{x})$ on a parameter $\vec{x}\in \mathbb{R}^n$ can be updated with the data $\vec{b}$ using the likelihood function $\pi(\vec{b}|\vec{x})$, more precisely, $\pi(\vec{x}|\vec{b})\,\propto\, \pi(\vec{b}|\vec{x})\pi(\vec{x})$. In many applications, the variable of interest $\vec{x}$ must satisfy some constraints: the pixels in images are bounded, the attenuation coefficients in a CT scan are nonnegative and the mass of an object might be bounded from above. This kind of information is often incorporated in the prior distribution $\pi(\vec{x})$, with the goal of making the posterior more accurate and explainable. 

Besides choosing common prior distributions that are naturally restricted to the constraint space, there are multiple other ways to add the constraints to the prior information. One of them is to truncate the prior distribution \cite{geweke1996bayesian}, where an unconstrained prior $\pi(\vec{x})$ is replaced by a prior proportional to $\pi(\vec{x})\mathbf{1}_{\set{C}}(\vec{x})$, where $\mathbf{1}_{\set{C}}(\vec{x})$ is $1$ if $\vec{x} \in \set{C}$ and $0$ otherwise. This is equivalent to truncating the posterior, resulting in a posterior proportional to $\pi(\vec{x}|\vec{b})\mathbf{1}_{\set{C}}(\vec{x})$. Samples from such a truncated posterior can be obtained using MCMC methods. Alternatively, the variable can be reparameterized to an unconstrained space. For example, for a positive variable $x$, we can write $x = e^{z}$ and define an unconstrained prior on $z$. If $z$ is normally distributed, then $x$ is referred to as the log-normal distribution \cite{bardsley2018computational}. These methods generally focus on the interior of the constraint set, while in many applications, the signals of interest lie on the boundary of the constraint set. For example, an image with at least one zero-valued pixel already lies on the boundary of the nonnegative orthant and in many applications, a lot of pixels are expected to be zero. However, if the distribution is described by a density, then the probability of having at least one zero-valued pixel is zero.

Here, we consider a method that focuses on the boundary of the constraint set by projecting samples from an unconstrained posterior onto the constraint set. Such methods have been analysed before in a Bayesian decision-theoretic framework \cite{sen2018constrained}. Furthermore, in the case of Bayesian linear inverse problems with Gaussian likelihood and prior with nonnegativity constraints, it has been noted \cite{bardsley2012mcmc, bardsley2020mcmc} that samples from a projected posterior can be obtained by sampling from constrained versions of randomized linear least squares problems. More precisely, let $\vec{x}$ be an unknown signal which is observed through a linear forward operator $A : \mathbb{R}^n \rightarrow \mathbb{R}^m$. Examples of linear forward operators include convolutions \cite{kaipio2006statistical} and CT-scans \cite{bangsgaard2021statistical}. The observations are inaccurate measurements of the form $\vec{b} = A\vec{x} + \vec{e}$ with error $\vec{e} \in \mathbb{R}^m$. The components of the error $\vec{e}$ are often modelled as independent and identically distributed Gaussian random variables, which results in a likelihood of the form $\pi(\vec{b}|\vec{x})\,\propto\, \exp\left(-\frac{\lambda}{2}\|A\vec{x}-\vec{b}\|_2^2\right)$. If we model our initial knowledge of $\vec{x}$ as the prior distribution $\pi(\vec{x})\,\propto\, \exp\left(-\frac{\delta}{2}\|L\vec{x}\|_2^2\right)$, then the posterior distribution satisfies
\begin{align} \label{eq:posterior_introduction}
    \pi(\vec{x}|\vec{b})\,\propto\, \pi(\vec{b}|\vec{x})\pi(\vec{x})\,\propto\, \exp\left(-\frac{\lambda}{2}\|A\vec{x}-\vec{b}\|_2^2 -\frac{\delta}{2}\|L\vec{x}\|_2^2\right).
\end{align}
In \cite{bardsley2020mcmc}, they observed that projecting \eqref{eq:posterior_introduction} onto the nonnegative orthant with respect to the norm $\|\cdot\|_{\lambda A^TA + \delta L^TL}$ is equivalent to solving the randomized constrained least squares problem
\begin{equation*}\label{eq:introduction}
\vec{x}^\star = \argmin_{\vec{x} \in \mathbb{R}_{+}^n}\left\{ \frac{\lambda}{2}\|A\vec{x}-\hat{\vec{b}}\|_2^2 + \frac{\delta}{2}\|L\vec{x}-\hat{\vec{c}}\|_2^2\right\},
\end{equation*}
where $\hat{\vec{b}} \sim \mathcal{N}(\vec{b}, \lambda^{-1}I)$ and $\hat{\vec{c}} \sim \mathcal{N}(\vec{0}, \delta^{-1}I)$.

In this paper, we generalize the approach from \cite{bardsley2012mcmc, bardsley2020mcmc} to general constraints, with a focus on polyhedral sets, and analyze the projected Gaussian posterior obtained through solving randomized constrained least squares problems. We derive a characterization for the projected distribution and show how this theory can be applied to Bayesian linear inverse problems. We also discuss how the projected posterior relates to a constrained prior and derive a Gibbs sampler for when the constraints are polyhedral cones. Finally, we use numerical experiments to look into the effect of the projection to Bayesian linear inverse problems. 

This paper is organized as follows. In Section \ref{sec_proj_Gaussian}, we discuss the theory behind the oblique projection of a multivariate Gaussian distribution onto a closed and convex set. In Section \ref{sec_lin_inverse_probs}, we discuss how to apply this projection framework to uncertainty quantification for linear inverse problems. Finally, in Section \ref{sec_num_examples}, we discuss the projection framework applied to deblurring and computer tomography.

\section{Projected multivariate Gaussian distribution}\label{sec_proj_Gaussian}
In this section, we describe a framework of projecting Gaussian distributions onto closed convex sets, where the projection is with respect to the norm induced by the precision matrix of the Gaussian. We describe theort for polyhedral constraint sets and give some explicit descriptions of the projected Gaussian for a few simple constraint sets.

\subsection{Oblique projection of Gaussians}

Throughout this section, unless otherwise stated, we make the following assumption.
\begin{assumption}\label{as_gauss}
The random vector $\vec{x}^\star \in \reals^n$ follows a Gaussian distribution with mean $\bm{\mu} \in \reals^n$ and covariance matrix $\Sigma \in \symm_{++}^{n}$, where $\symm_{++}^{n}$ denotes the set of $n \times n$ symmetric positive definite matrices.
\end{assumption}

For $\vec{x}^\star$ as defined by Assumption \ref{as_gauss}, we can equivalently write $\vec{x}^\star = \bm{\mu} + \Sigma \hat{\vec{w}}$, where $\hat{\vec{w}} \sim \mathcal{N}(\vec{0},\Sigma^{-1})$. Define the quadratic function $g(\vec{x}; \vec{w}) := \frac12 \vec{x}^T\Sigma^{-1}\vec{x} - \vec{x}^T(\Sigma^{-1}\bm{\mu} + \vec{w})$, then $\vec{x}^\star$ satisfies $\nabla_{\vec{x}} g(\vec{x}^\star; \hat{\vec{w}}) = \vec{0}$. Therefore, $\vec{x}^\star$ is the solution to a randomized quadratic optimization problem of the form
\begin{equation}\label{e_gauss_opt}
    \argmin_{\vec{x}\in \reals^n} \left\{\frac12 \vec{x}^T\Sigma^{-1}\vec{x} - \vec{x}^T(\Sigma^{-1}\bm{\mu} + \hat{\vec{w}})\right\}.
\end{equation}
Thus, we can sample from $\vec{x}^\star$ by repeatedly solving optimization problem \eqref{e_gauss_opt} for different samples of $\hat{\vec{w}}$.

Now consider constraining optimization problem \eqref{e_gauss_opt} to a closed convex set $\set{C} \subseteq \reals^n$. The resulting optimization problem has the form
\begin{equation*}
    \vec{z}^\star = \argmin_{\vec{z}\in \set{C}} \left\{\frac12 \vec{z}^T\Sigma^{-1}\vec{z} - \vec{z}^T(\Sigma^{-1}\bm{\mu} + \hat{\vec{w}})\right\},\quad \hat{\vec{w}} \sim \mathcal{N}(\vec{0}, \Sigma^{-1}),
\end{equation*}
or equivalently
\begin{equation*}
    \vec{z}^\star = \argmin_{\vec{z}\in \set{C}} \frac12\|\vec{z} - (\bm{\mu} + \Sigma \hat{\vec{w}})\|^2_{\Sigma^{-1}} = \argmin_{\vec{z}\in \set{C}} \frac12\|\vec{z} - \vec{x}^\star\|^2_{\Sigma^{-1}},
\end{equation*}
i.e., $\vec{z}^\star$ is the oblique projection of $\vec{x}^\star \sim \mathcal{N}(\bm{\mu}, \Sigma)$ onto $\set{C}$ with respect to the norm induced by the precision matrix $\Sigma^{-1}$. This obliquely projected Gaussian distribution is the main object of study in this work and is summarized in the following definition.

\begin{definition}
Under Assumption \ref{as_gauss}, the oblique projection of $\mathcal{N}(\bm{\mu}, \Sigma)$ onto a closed, convex set $\set{C} \subset \reals^n$ is the oblique projection of $\vec{x}^\star$ onto the set $\set{C}$ with respect to the norm induced by the precision matrix $\Sigma^{-1}$, that is,
\begin{equation}\label{def_projection}
\vec{z}^\star = \Pi_{\set{C}}^{\Sigma^{-1}}(\vec{x}^\star) := \argmin_{\vec{z}\in \set{C}} \frac12\|\vec{x}^\star - \vec{z}\|^2_{\Sigma^{-1}}.
\end{equation}

\end{definition}

Because the precision matrix $\Sigma^{-1}$ is positive definite, the oblique projection onto a closed, convex set $\Pi_{\set{C}}^{\Sigma^{-1}}$ is well-defined and continuous, hence measurable. Therefore, the random vector \eqref{def_projection} is well defined with distribution
\begin{equation}\label{e_projected_distribution}
    \mathbb{P}(\Pi_{\set{C}}^{\Sigma^{-1}}(\vec{x}^\star) \in \set{E}) = \mathbb{P}\left(\vec{x}^\star \in \left[\Pi_{\set{C}}^{\Sigma^{-1}}\right]^{-1}(\set{E})\right),
\end{equation}
for a measurable set $\set{E}$ and where $\left[\Pi_{\set{C}}^{\Sigma^{-1}}\right]^{-1}$ denotes the inverse image of $\Pi_{\set{C}}^{\Sigma^{-1}}$.

Solving the optimization problem \eqref{def_projection} can be computationally expensive, however, the inverse mapping can be analyzed more easily. The optimality condition of the oblique projection is given by
\begin{equation*}\label{e_opt_cond}
\vec{0} \in \Sigma^{-1}(\vec{z}^\star-\vec{x}^\star) + N_{\set{C}}(\vec{z}^\star),
\end{equation*}
or equivalently
\begin{equation*}
\vec{x}^\star \in \vec{z}^\star + \Sigma N_{\set{C}}(\vec{z}^\star),
\end{equation*}
where $N_{\set{C}}(\vec{z}) = \partial I_{\set{C}}(\vec{z}) = \{\vec{v}\,|\, \vec{v}^T(\vec{y}-\vec{z}) \leq 0, \forall \vec{y} \in \set{C}\}$ is the normal cone associated with $\set{C}$.

From optimality condition \eqref{e_opt_cond} and equation \eqref{e_projected_distribution} we get that the distribution of the projected Gaussian can be described by
\begin{equation}\label{eq_inversion}
    \mathbb{P}\left(\Pi_{\set{C}}^{\Sigma^{-1}}(\vec{x}^\star) \in \set{E}\right) = \mathbb{P}\left(\vec{x}^\star \in \bigcup_{\vec{z}\in \set{E}}\left[\vec{z} + \Sigma N_{\set{C}}(\vec{z})\right]\right),
\end{equation}
where $\set{E}$ is a measurable set. 

Because of the projection, a lot of the mass gets projected onto the boundary of the constraint set. If $\set{C}$ is a polyhedral set, then the probability on the boundary can be described as in the following lemma.

\begin{lemma}\label{lemma_densities}
Under Assumption \ref{as_gauss}, if $\set{C}$ is a polyhedral set and $\set{F}$ is a face of $\set{C}$, then for any measurable set $\set{E} \subseteq \relint(\set{F})$, where $\relint(\set{F})$ denotes the relative interior of the face $\set{F}$, we have
\begin{equation*}
    \mathbb{P}\left(\Pi_{\set{C}}^{\Sigma^{-1}}(\vec{x}^\star) \in \set{E}\right) = \int_{\set{E}} \pi_{\set{F}}(\vec{z}) \text{\normalfont d}\vec{z},
\end{equation*}
where $\pi_{\set{F}}(\vec{z})$ is the $\dim(\set{F})$-dimensional density
\begin{equation*}
    \pi_{\set{F}}(\vec{z}) = \int_{\Sigma N_\set{C}(\vec{z})} \pi_{\vec{x}^\star}(\vec{z}+\vec{v}) \text{\normalfont  d}\vec{v}.
\end{equation*}
\end{lemma}
\begin{proof}
For any $\vec{z} \in \relint(\set{F})$, the normal cone $N_{\set{C}}(\vec{z})$ is independent of $\vec{z}$, therefore we can write \eqref{eq_inversion} as
\begin{align*}
\mathbb{P}\left(\Pi_{\set{C}}^{\Sigma^{-1}}(\vec{x}^\star) \in \set{E}\right)
&= \int_{\set{E} + \Sigma N_{\set{C}}} \pi_{\vec{x}^{\star}}(\vec{w})\text{d}\vec{w}\\
&= \int_{\set{E}}\int_{\Sigma N_{\set{C}}}\pi_{\vec{x}^\star}(\vec{z}+\vec{v})\text{d}\vec{v}\text{d} \vec{z},
\end{align*}
where the final decomposition is valid due to Lemma \ref{lem:basis_transform} in Appendix \ref{subsec:appendix_proof}.
\end{proof}

Lemma \ref{lemma_densities} above states that the projected Gaussian distribution onto a polyhedral set consists of a mixture of various densities of different dimensions on all of the faces of the polyhedral set.

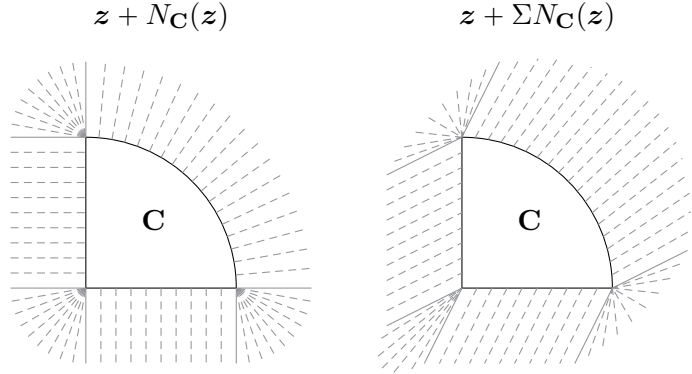
\begin{figure}
    \centering
    \begin{tikzpicture}[scale=2]
        \draw (1,0) arc (0:90:1);
        \draw[black] (0,0) -- (0,1);
        \draw[black] (0,0) -- (1,0);
        
        \draw[gray] (0,0) -- (-0.5,0);
        \draw[gray] (0,0) -- (0,-0.5);
        \draw[gray] (1,0) -- (1.5,0);
        \draw[gray] (1,0) -- (1,-0.5);
        \draw[gray] (0,1) -- (0,1.5);
        \draw[gray] (0,1) -- (-0.5,1);
        
        \foreach \y in {0.1,0.2,...,0.9}{\draw[gray, densely dashed] (-0.5,\y) -- (0,\y);}
        \foreach \x in {0.1,0.2,...,0.9}{\draw[gray, densely dashed] (\x,0) -- (\x,-0.5);}
        \foreach \phi in {5,10,...,85}{\draw[gray, densely dashed] ({cos(\phi)},{sin(\phi)}) -- ({1.5*cos(\phi)},{1.5*sin(\phi)});}
        \foreach \phi in {10,20,...,80}{\draw[gray, densely dashed] ({-0.5*cos(\phi)},{-0.5*sin(\phi)}) -- (0,0);}
        \foreach \phi in {10,20,...,80}{\draw[gray, densely dashed] ({-0.5*cos(\phi)},{1+0.5*sin(\phi)}) -- (0,1);}
        \foreach \phi in {10,20,...,80}{\draw[gray, densely dashed] ({1+0.5*cos(\phi)},{-0.5*sin(\phi)}) -- (1,0);}
        
        \node[] at (0.5,1.8) {$\vec{z} + N_{\set{C}}(\vec{z})$};
        \node[] at (0.45,0.45) {$\set{C}$};
        \node[] at (3,1.8) {$\vec{z} + \Sigma N_{\set{C}}(\vec{z})$};
        \node[] at (2.95,0.45) {$\set{C}$};
        
        \draw (3.5,0) arc (0:90:1);
        \draw[black] (2.5,0) -- (2.5,1);
        \draw[black] (2.5,0) -- (3.5,0);
        
        \draw[gray] (2.5,0) -- (2,-0.25);
        \draw[gray] (2.5,0) -- (2.25,-0.5);
        \draw[gray] (3.5,0) -- (4,0.25);
        \draw[gray] (3.5,0) -- (3.25,-0.5);
        \draw[gray] (2.5,1) -- (2.75,1.5);
        \draw[gray] (2.5,1) -- (2,0.75);
        
        \foreach \y in {0.1,0.2,...,0.9}{\draw[gray, densely dashed] (2,{\y-0.25}) -- (2.5,\y);}
        \foreach \x in {0.1,0.2,...,0.9}{\draw[gray, densely dashed] ({\x+2.5},0) -- ({\x+2.25},-0.5);}

        \foreach \phi in {5,10,...,85}{\draw[gray, densely dashed] ({2.5+cos(\phi)},{sin(\phi)}) -- ({2.5+1.5*cos(\phi)+0.25*sin(\phi)},{1.5*sin(\phi) + 0.25*cos(\phi)});}
        
        \foreach \phi in {10,27.5,...,80}{\draw[gray, densely dashed] ({2.5-0.5*cos(\phi)-0.25*sin(\phi)},{-0.5*sin(\phi)-0.25*cos(\phi)}) -- (2.5,0);}
        \foreach \phi in {10,20,...,80}{\draw[gray, densely dashed] ({2.5-0.5*cos(\phi)+0.25*sin(\phi)},{1+0.5*sin(\phi)-0.25*cos(\phi)}) -- (2.5,1);}
        \foreach \phi in {10,20,...,80}{\draw[gray, densely dashed] ({2.5+1+0.5*cos(\phi)-0.25*sin(\phi)},{-0.5*sin(\phi)+0.25*cos(\phi)}) -- (3.5,0);}
        
    \end{tikzpicture}

    \caption{Visualization of obliquely projecting a density onto the different boundary components of a quarter disc.}
    \label{fig-projection}
\end{figure}

Extending Lemma \ref{lemma_densities} to general closed, convex sets $\set{C}$ is a more complicated procedure, but the idea of densities of different dimensions on different parts of the constraint set is the same. For example, consider the quarter disc as illustrated in Figure \ref{fig-projection}. The projection has no effect on the interior of the domain, therefore, there is a two-dimensional density on the interior of the set. The normal cones at the three corners are two dimensional, therefore, positive mass gets projected onto each of the corners. The normal cone at the straight and curved lines are one-dimensional, resulting in one-dimensional densities on these parts. Figure \ref{fig:quarter_disc}, shows these densities, except for the corner points. To further illustrate this, we will consider a few examples for which we can compute analytical expression for the densities.

\subsection{Analytic examples}\label{subsec:analytical_examples}
Computing analytical expressions for the densities on parts of the surface of the constraint set is in general intractable, but can be computed when the dimension of the normal cone $N_{\set{C}}(\vec{x})$ is at most one-dimensional. In this subsection, we consider two such cases, half-space and disc constraints, and give analytical expressions for their densities.

\subsubsection{Halfspace}
Suppose that $\set{C}$ is a halfspace defined by $\set{C} = \{\vec{x}\in \reals^n \,|\, \vec{a}^T\vec{x} \leq b\}$, where $\vec{a} \in \reals^n$ is a nonzero normal vector and $b \in \reals$. Denote by $F \in \reals^{n\times(n-1)}$ a matrix whose columns form an orthonormal basis for $\nullspace(\vec{a}^T)$ and let $\vec{x}_0 \in \mathbb{R}^n$ satisfy $\vec{a}^T\vec{x}_0 = b$, then the boundary of the halfspace $\set{C}$ can be parameterized by $\vec{x}_0 + F\vec{u}$ for $\vec{u}\in \reals^{n-1}$. Under Assumption \ref{as_gauss},  the $n-1$ dimensional density of $\Pi_{\set{C}}^{\Sigma^{-1}}(\vec{x}^\star)$ on the boundary of $\set{C}$ is given by
\begin{equation}\label{e-halfspace}
\pi_{\bd}(\vec{u}) = \pi_{\vec{x}^{\star}}(\vec{x}_0+F\vec{u})\sqrt{\frac{\vec{a}^T\Sigma \vec{a}}{\vec{a}^T\vec{a}}} \exp\left(\gamma^2\right)\sqrt{\frac{\pi}{2}} \erfc\left(\gamma \right),
\end{equation}
where $\gamma = \frac{2(b-\vec{a}^T\bm{\mu})}{\sqrt{8\vec{a}^T\Sigma \vec{a}}}$ and $\erfc$ is the complementary error function.

\subsubsection{Unit disc}
Suppose that $\set{C}$ is a unit disc defined by $\set{C} = \{\vec{x} \in \reals^2 \,|\, \|\vec{x}\|_2 \leq 1\}$ with boundary parameterization $\vec{n}(u) := (\cos(u), \sin(u))^T$ for $u \in [0, 2\pi)$. The density on the boundary of the disc is given by
\begin{equation}\label{e-disc}
\pi_{\bd}(u) = \pi_{\vec{x}^{\star}}(\vec{n}(u))\left(\frac{\det(\Sigma)}{\alpha} + \frac{\sqrt{\pi}}{\sqrt{8}a^{3/2}}(2\alpha K(u) - \beta \det(\Sigma))\exp\left(\frac{\beta^2}{8\alpha} \right)\erfc\left(\frac{\beta}{\sqrt{8\alpha}}\right)\right),
\end{equation}
where $\alpha = \vec{n}(u)^T\Sigma \vec{n}(u)$, $\beta = 2\vec{n}(u)^T(\vec{n}(u)-\bm{\mu})$, $K(u) = \det [\Sigma \vec{n}(u)\ R\vec{n}(u)]$ and $R = \begin{bmatrix}0 & -1 \\ 1 & 0\end{bmatrix}$.

\begin{figure}
    \centering
    \includegraphics[width = 0.5\textwidth]{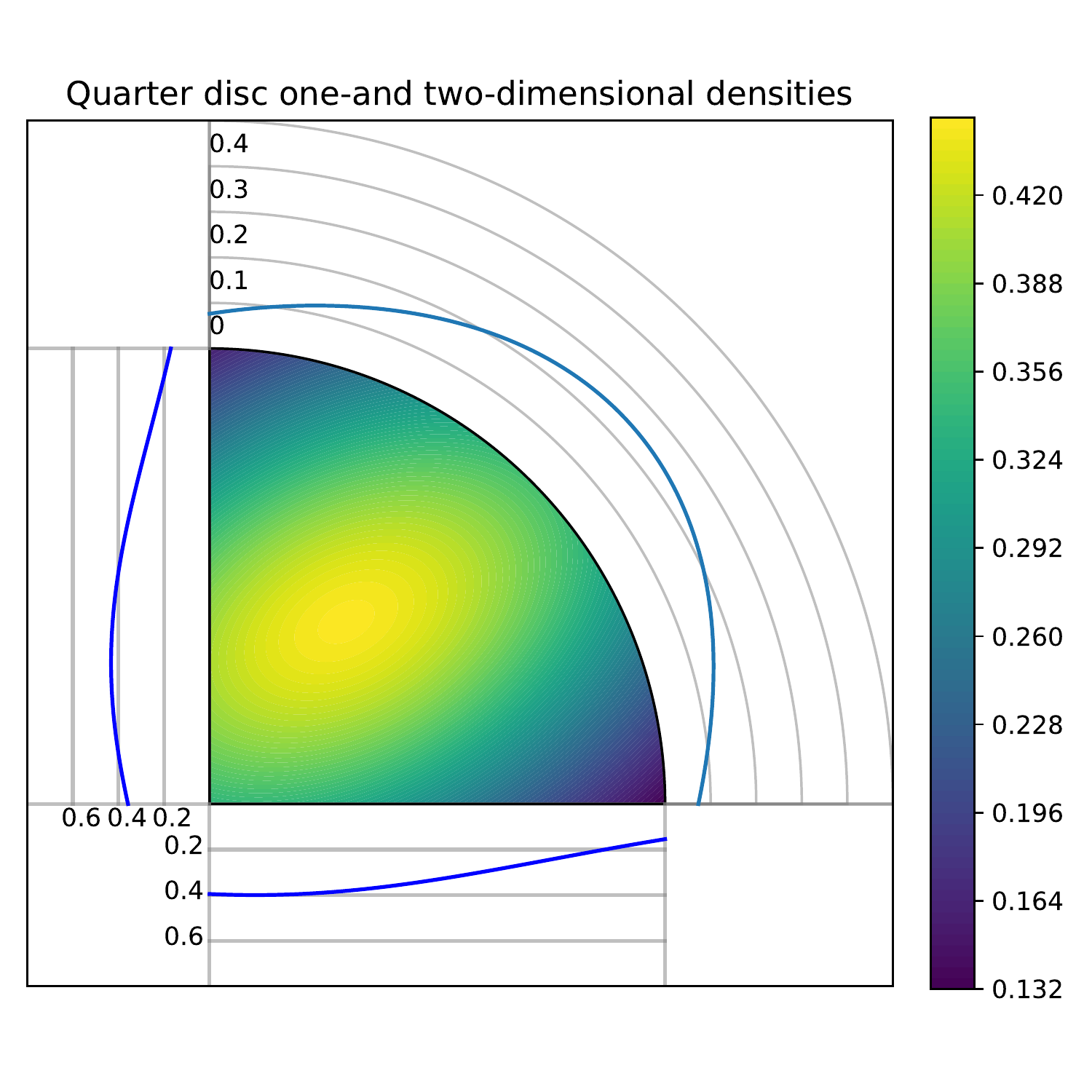}
    \caption{The one-and two-dimensional densities of the oblique projection of a Gaussian onto a quarter disc. The zero-dimensional densities in the corners are not shown.}
    \label{fig:quarter_disc}
\end{figure}

These computations are all based on being able to compute the integral of a Gaussian over a one-dimensional normal cone and the computations can be found in Appendix \ref{subsec:appendix_analytical_examples}. Furthermore, these distributions can be mixed. Figure \ref{fig:quarter_disc} shows the one-and two-dimensional densities on the boundary of a quarter disc as described by the Equations \eqref{e-halfspace} and \eqref{e-disc}.

\subsection{Boundary properties}
As seen in the analytical examples in Subsection \ref{subsec:analytical_examples}, a lot of mass is projected onto the boundary of the constraint set $\set{C}$. Nevertheless, there will always be a positive probability on the relative interior of the constraint set. The following lemma shows that both the boundary and relative interior of the constraint set will always have positive probability. Although the lemma is formulated in terms of the oblique projection of a Gaussian distribution with respect to its precision matrix, the proof extends to the projection with respect to any positive definite matrix and any continuous distribution whose support is $\reals^n$.
\begin{lemma}\label{lem:boundary_interior_positive_probability}
Under Assumption \ref{as_gauss}, if $\set{C} \subset \reals^n$ is a non-empty, closed, convex set, then the probability of being on the boundary of $\set{C}$ is positive and the probability of being in the relative interior of $\set{C}$ is positive, i.e., 
\begin{equation*}
\mathbb{P}(\Pi^{\Sigma^{-1}}_{\set{C}}(\vec{x}^\star) \in \bd(\set{C})) > 0\quad \text{and}\quad \mathbb{P}(\Pi^{\Sigma^{-1}}_{\set{C}}(\vec{x}^\star) \in \relint(\set{C})) > 0.
\end{equation*}
\end{lemma}
\begin{proof}
    Because $\set{C} \subset \mathbb{R}^n$ and closed, $\mathbb{R}^n \backslash \set{C}$ is non-empty and open. Therefore 
    \begin{equation*}
    \mathbb{P}(\Pi^{\Sigma^{-1}}_{\set{C}}(\vec{x}^\star) \in \bd(\set{C})) \geq \mathbb{P}(\vec{x}^\star \in \mathbb{R}^n \backslash \set{C}) > 0.
    \end{equation*}
    Let $\vec{z} \in \relint(\set{C})$, then there exists $\epsilon > 0$ such that $B_\epsilon(\vec{z}) \cap \aff(\set{C}) \subseteq \set{C}$, where $B_\epsilon(\vec{z})$ is the closed ball of radius $\epsilon$ around $\vec{z}$ and $\aff(\set{C})$ is the affine hull of $\set{C}$. By Lemma \ref{lem:basis_transform}, $B_\epsilon(\vec{z}) \cap \aff(\set{C}) + \Sigma \aff(\set{C})^\bot$ is an $n$-dimensional convex set. Therefore we can conclude that
    \begin{equation*}
    \mathbb{P}(\Pi^{\Sigma^{-1}}_{\set{C}}(\vec{x}^\star) \in \relint(\set{C})) \geq \mathbb{P}(\vec{x}^\star \in B_\epsilon(\vec{z}) \cap \aff(\set{C}) + \Sigma \aff(\set{C})^\bot)> 0.
    \end{equation*}
\end{proof}

Typically, one would like to compute a single point estimate from the posterior and analyze the uncertainty of that estimate. A common choice for point estimate is the MAP (Maximum a Posteriori), which is the point of maximum posterior density. However, due to the projected posterior consisting of a mixture of different dimensional densities, the MAP estimate is undefined. Another common choice for point estimate is the mean of the posterior, however, as the following theorem shows, the mean of the projected posterior lies inside the relative interior of the constraint set. It is therefore not always a suitable point estimate if one is interested in the boundary of the constraint set. This motivates our choice of the componentwise median as point estimate in the numerical examples in Section \ref{sec_num_examples}.
\begin{theorem}
Under Assumption \ref{as_gauss}, if $\set{C} \subset \mathbb{R}^n$ is a non-empty, closed, convex set, then
\begin{equation*}
\mathbb{E}[\Pi^{\Sigma^{-1}}_{\set{C}}(\vec{x}^\star)] \in \relint(\set{C}).
\end{equation*}
\end{theorem}
\begin{proof}
By Lemma \ref{lem:boundary_interior_positive_probability}, $p := \mathbb{P}(\Pi^{\Sigma^{-1}}_{\set{C}}(\vec{x}^\star) \in \relint(\set{C})) > 0$. Let $\vec{z}^\star = \Pi^{\Sigma^{-1}}_{\set{C}}(\vec{x}^\star)$, then consider the decomposition of the expectation on the relative interior, i.e.,
\begin{equation*}
\mathbb{E}[\vec{z}^\star] = p\mathbb{E}[\vec{z}^\star|\vec{z}^\star \in \relint(\set{C})] + (1-p)\mathbb{E}[\vec{z}^\star| \vec{z}^\star \in \relbd(\set{C})].
\end{equation*}
We have that $\mathbb{E}[\vec{z}^\star|\vec{z}^\star \in \relint(\set{C})] \in \relint(\set{C})$ by the convexity of the relative interior and $\mathbb{E}[\vec{z}^\star| \vec{z}^\star \in \relbd(\set{C})] \in \set{C}$. Combined with the fact \cite[Lemma 2.1.6]{hiriart2004fundamentals} that for $p \in (0, 1]$, $p\relint(\set{C}) + (1-p)\set{C} \subseteq \relint(\set{C})$, we can conclude that $\mathbb{E}[\vec{z}^\star] \in \relint(\set{C})$.
\end{proof}

\subsection{Gaussian decomposition}
The two analytical examples in Subsection \ref{subsec:analytical_examples} illustrate that the densities on the boundary of the constraint set can become quite complicated, even for relatively simple examples. However, note that the density on the boundary of a half-space in \eqref{e-halfspace} is proportional to the original Gaussian density, while on the boundary of a disc in \eqref{e-disc}, the additional factor depends on the surface coordinate. The property that the densities are proportional to the original Gaussian as observed in the half-space example can be generalized to polyhedral sets as stated in the following theorem, which is a key observation for the derivation of the Gibbs sampler in Subsections \ref{ss_constrained_prior} and \ref{ss_Gibbs}.

\begin{theorem}\label{thm_gaussian_decomposition}
Under Assumption \ref{as_gauss}, if $\set{C} \subseteq \reals^n$ is a polyhedral set, then the density of the projected Gaussian $\Pi_{\set{C}}^{\Sigma^{-1}}(\vec{x}^\star)$ on the relative interior of any face of $\set{C}$ is proportional to the density of $\vec{x}^\star$.
\end{theorem}
\begin{proof}

    Consider a face $\set{F}$ of the polyhedral set $\set{C} \subseteq \reals^n$ with linearly independent normal vectors $\vec{a}_1, \dots, \vec{a}_k \in \reals^n$, where $k$ is the dimension of the face $\set{F}$, and let $F \in \reals^{n \times (n-k)}$ be a matrix whose columns form a basis for the space orthogonal to the span of the normal vectors $\vec{a}_1, \dots, \vec{a}_k$, i.e., $F^T\vec{a}_i = \vec{0}$ for all $i = 1,\dots, k$. If $\vec{x}_0$ is any point on $\set{F}$, then each point $\vec{z} \in \set{F}$ can be written as $\vec{x}_0 + F\vec{u}$ for some $\vec{u}$ and the corresponding normal cone can be parameterized as $\sum_{i=1}^{k}t_i\vec{a}_i$ for $t_i \geq 0$ for all $i = 1,\dots,k$. The density $\pi_{\set{F}}$ on the face then satisfies
    \begin{equation*}
    	\pi_{\set{F}}(\vec{u}) \,\propto\, \int_{0}^{\infty}\cdots\int_{0}^{\infty}\pi_{\vec{x}^{\star}}\left(\vec{x}_0 + F\vec{u} + \sum_{i=1}^{k}t_i\Sigma \vec{a}_i\right) \text{d}t_1\cdots \text{d}t_k.
    \end{equation*}
    Note that
    \begin{align*}
    	-&2\log\left(\frac{\pi_{\vec{x}^{\star}}\left(\vec{x}_0 + F\vec{u} + \sum_{i=1}^{k}t_i\Sigma \vec{a}_i\right)}{\pi_{\vec{x}^{\star}}\left(\vec{x}_0 + F\vec{u
    	}\right)}\right) \\ &=
    	 2 \left(\sum_{i=1}^{k}t_i \vec{a}_i\right)^T\left(\vec{x}_0 + F\vec{u}-\bm{\mu}\right) + \left(\sum_{i=1}^{k}t_i\Sigma \vec{a}_i\right)^T\Sigma^{-1}\left(\sum_{i=1}^{k}t_i\Sigma \vec{a}_i\right).
    \end{align*}
    The only term that depends on $\vec{u}$ vanishes because
    \begin{equation*}
    	\left(\sum_{i=1}^{k}t_i \vec{a}_i\right)^TF\vec{u} = \sum_{i=1}^{k}t_i \vec{a}_i^TF\vec{u} = 0,
    \end{equation*}
    and hence
    \begin{equation*}
    	\pi_{\vec{x}^{\star}}\left(\vec{x}_0 + F\vec{u} + \sum_{i=1}^{k}t_i\Sigma \vec{a}_i\right) \,\propto\, c(t_1,\dots t_k)\pi_{\vec{x}^{\star}}\left(\vec{x}_0 + F\vec{u}\right),
    \end{equation*}
    where $c(t_1,\dots t_k)$ only depends on the coefficients $t_i$. It follows that
    \begin{equation*}\label{e_aff_gauss}
    	\pi_{\set{F}}(\vec{u})\, \propto\, \int_{0}^{\infty}\dots\int_{0}^{\infty}\pi_{\vec{x}^{\star}}\left(\vec{x}_0 + F\vec{u} + \sum_{i=1}^{k}t_i\Sigma \vec{a}_i\right) \text{d}t_1\dots \text{d}t_k\, \propto\, \pi_{\vec{x}^{\star}}\left(\vec{x}_0 + F\vec{u}\right),
    \end{equation*}
    i.e., the density of the projected Gaussian on a face of the polyhedral set is proportional to the unprojected Gaussian.
\end{proof}

\section{Bayesian linear inverse problems with constraints}\label{sec_lin_inverse_probs}
In this section, we describe how to apply the theory for projected Gaussian distributions in Section \ref{sec_proj_Gaussian} to linear inverse problems. We discuss the randomized optimization problems to be solved to obtain samples from the projected posteriors and how the projected posterior relates to a constrained prior. Finally, we define a Bayesian hierarchical model for the linear inverse problem and derive a Gibbs sampler for that model in the special case where the constraint set $\set{C}$ is a polyhedral cone.

\subsection{Bayesian model}
Let us now consider the problem of recovering a signal $\vec{x} \in \reals^n$ from noisy observations $\vec{b} = A\vec{x} + \vec{e}$ using a linear forward operator $A : \mathbb{R}^n \rightarrow \mathbb{R}^m$ and with noise $\vec{e} \sim \mathcal{N}(\vec{0}, \Sigma_{\vec{e}})$ for $\Sigma_{\vec{e}}\in  \symm_{++}^n$. This results in the likelihood function
\begin{equation*}\label{e_likelihood}
    \pi(\vec{b}\,|\,\vec{x})\,\propto\, \exp\left(-\frac12\|A\vec{x}-\vec{b}\|^2_{\Sigma_{\vec{e}}^{-1}}\right).
\end{equation*}
Furthermore, assume a priori that $\vec{x} \sim \mathcal{N}(\vec{0}, \Sigma_{\vec{x}})$ with $\Sigma_{\vec{x}} \in \symm_{++}^n$, i.e.,
\begin{equation*}
    \pi(\vec{x})\,\propto\, \exp\left(-\frac12\|\vec{x}\|^2_{\Sigma_{\vec{x}}^{-1}}\right).
\end{equation*}
From Bayes' theorem we obtain the posterior distribution
\begin{equation}\label{e_posterior_density}
    \pi(\vec{x}\,|\,\vec{b}) \,\propto\, \pi(\vec{b}\,|\,\vec{x})\pi(\vec{x}) \,\propto\, \exp\left(-\frac12\|A\vec{x}-\vec{b}\|^2_{\Sigma_{\vec{e}}^{-1}} -\frac12\|\vec{x}\|^2_{\Sigma_{\vec{x}}^{-1}}\right),
\end{equation}
which is a Gaussian distribution with  covariance $\Sigma_{\vec{x}|\vec{b}} = (\Sigma_{\vec{x}}^{-1} + A^T\Sigma_{\vec{e}}^{-1}A)^{-1}$ and mean $\bm{\mu} = \Sigma_{\vec{x}|\vec{b}}A^T\Sigma_{\vec{e}}^{-1}\vec{b}$. Now consider a closed, convex set $\set{C} \subset \reals^n$. To sample from the projected posterior, instead of solving
\begin{equation*}
    \argmin_{\vec{x}\in \set{C}} \left\{\frac12 \vec{x}^T(\Sigma_{\vec{x}}^{-1} + A^T\Sigma_{\vec{e}}^{-1}A)\vec{x} - \vec{x}^T(A^T\Sigma_{\vec{e}}^{-1}\vec{b} + \hat{\vec{w}})\right\},
\end{equation*}
where $\hat{\vec{w}} \sim \mathcal{N}(\vec{0}, \Sigma_{\vec{x}|\vec{b}}^{-1})$, we can equivalently solve
\begin{equation}\label{e_randomized_constrained_least_squares}
    \argmin_{\vec{x}\in \set{C}} \left\{\frac12 \|A\vec{x}-\hat{\vec{b}}\|^2_{\Sigma_{\vec{e}}^{-1}} + \frac12 \|\vec{x}-\hat{\vec{c}}\|_{\Sigma_{\vec{x}}^{-1}}^2\right\},
\end{equation}
where $\hat{\vec{b}} \sim \mathcal{N}(\vec{b}, \Sigma_{\vec{e}})$ and $\hat{\vec{c}} \sim \mathcal{N}(\vec{0}, \Sigma_{\vec{x}})$, i.e., solving randomized constrained linear least squares problems or randomized MAP estimates.

For a simple case where $\vec{e} \sim \mathcal{N}(\vec{0}, (\lambda I)^{-1})$ and $\vec{x} \sim \mathcal{N}(\vec{0}, (\delta L^T L)^{-1})$ with hyperparameters $\lambda,\delta > 0$ and $L$ is a full-rank matrix, the posterior \eqref{e_posterior_density} simplifies to
\begin{equation}\label{e_simplified_posterior}
    \pi(\vec{x}\,|\,\vec{b}) \,\propto\, \pi(\vec{b}\,|\,\vec{x})\pi(\vec{x}) \,\propto\, \exp\left(-\frac{\lambda}{2}\|A\vec{x}-\vec{b}\|^2_{2} -\frac{\delta}{2} \|L\vec{x}\|_{2}^2\right).
\end{equation}
Optimization problem \eqref{e_randomized_constrained_least_squares} then simplifies to
\begin{equation}\label{e_randomized_constrained_least_squares_simple}
    \argmin_{\vec{x}\in \set{C}} \left\{ \frac{\lambda}{2} \|A\vec{x}-\hat{\vec{b}}\|^2_{2} + \frac{\delta}{2} \|L\vec{x}-\hat{\vec{c}}\|_{2}^2 \right\},
\end{equation}
where $\hat{\vec{b}} \sim \mathcal{N}(\vec{b},(\lambda I)^{-1})$ and $\hat{\vec{c}} \sim \mathcal{N}(\vec{0}, (\delta I)^{-1})$.

Thus, by repeatedly solving the randomized constrained least squares problem \eqref{e_randomized_constrained_least_squares_simple}, we obtain samples from the obliquely projected posterior onto the set $\set{C}$.

\subsection{Constrained prior}\label{ss_constrained_prior}
Although we defined the posterior implicitly through modification of the unconstrained posterior, we can also interpret this modification as having implicitly defined a constrained prior. Let $\set{C} \subseteq \reals^n$ be a polyhedral set and define the projected posterior onto $\set{C}$ as in \eqref{def_projection}. By Theorem \ref{thm_gaussian_decomposition}, the projected posterior density on any face $\set{F}$ of $\set{C}$ is proportional to the unprojected posterior. Because constraining the signal is based on prior information, we will assume that the likelihood $\pi(\vec{b}|\vec{x})$ is not affected by the projection. Hence, by Bayes' formula, we obtain the prior density on the face,
\begin{align}\label{e:constrained_prior}
    \pi_{\vec{x},\set{F}}(\vec{u})&\,\propto\,  \frac{\pi_{\vec{x}|\vec{b},\set{F}}(\vec{u})}{\pi_{\vec{b}|\vec{x}}(\vec{x}_0 + F\vec{u})}\,\propto\, \frac{\exp\left(-\frac12\|A\vec{x}-\vec{b}\|^2_{\Sigma_{\vec{e}}^{-1}} -\frac12\|\vec{x}\|^2_{\Sigma_{\vec{x}}^{-1}}\right)}{\exp\left(-\frac12\|A\vec{x}-\vec{b}\|^2_{\Sigma_{\vec{e}}^{-1}}\right)}\\
    \nonumber&\,\propto\, \exp\left(-\frac12\|\vec{x}\|^2_{\Sigma_{\vec{x}}^{-1}}\right) \,\propto\,\pi_{\vec{x}}(\vec{x}_0 + F\vec{u}),
\end{align}
where we denote by $\vec{x} = \vec{x}_0 + F\vec{u}$ the point of the face $\set{F}$ parameterized by $\vec{u}$.

This shows that the corresponding constrained prior is proportional to the unconstrained prior on any face of $\set{C}$ and therefore has a similar structure on $\set{C}$ as the posterior.
However, this constrained prior is generally not the same as obliquely projecting the prior onto the constraint set.

\subsection{Bayesian hierarchical model and Gibbs sampler}\label{ss_Gibbs}
Let us now consider adding priors to the hyperparameters $\lambda$ and $\delta$. To exploit that Gaussian and Gamma distributions are conjugate, let the hyperpriors be $\lambda \sim \Gamma(\alpha_{\lambda}, \beta_{\lambda})$ and $\delta \sim \Gamma(\alpha_{\delta}, \beta_{\delta})$, i.e.,
\begin{align*}\label{e:hyperpriors}
    \pi(\lambda) &\,\propto\, \lambda^{\alpha_{\lambda}-1}\exp(-\beta_{\lambda}\lambda), \quad\text{ for } \lambda > 0 \quad \text{and}\\
    \nonumber
    \pi(\delta) &\,\propto\, \delta^{\alpha_{\delta}-1}\exp(-\beta_{\delta}\delta), \quad\ \text{ for } \delta > 0.
\end{align*}
Combining these hyperpriors with the likelihood and prior of the previous section results in a Bayesian hierarchical model. A common method for sampling for the signal $\vec{x}$ and hyperparameters $\lambda$ and $\delta$ is a hierarchical Gibbs sampler \cite[Algorithm 5.1]{bardsley2018computational}. For the unconstrained setting of \eqref{e_simplified_posterior}, it can be derived that
\begin{align*}
\pi_{x,\lambda,\delta | b}(x, \lambda, \delta) 
\,\propto\,&\ \lambda^{m/2+\alpha_{\lambda} - 1}\delta^{n/2+\alpha_{\delta} - 1}\\
       &\times\exp\left(-\frac{\lambda}{2}\|A\vec{x} - \vec{b}\|_2^2 - \frac{\delta}{2}\|L\vec{x}\|_2^2 - \beta_{\lambda}\lambda - \beta_{\delta}\delta\right),
\end{align*}
from which it follows that
\begin{align}
    \nonumber
    \lambda\, |\, \vec{x},\vec{b} &\sim \Gamma(m/2 + \alpha_{\lambda}, \frac{1}{2}\|A\vec{x}-\vec{b}\|_2^2 + \beta_\lambda),\\
    \delta\, |\, \vec{x},\vec{b} &\sim \Gamma(n/2 + \alpha_{\delta}, \frac{1}{2}\|L\vec{x}\|_2^2 + \beta_\delta),
    \label{eq:delta_unconstrained}
\end{align}
and $\vec{x}\,|\, \lambda, \delta, \vec{b}$ is a Gaussian described by \eqref{e_simplified_posterior}. A Gibbs sampler then alternates among sampling from these conditional distributions.

Now for the constrained setting, let $\set{C}$ be a polyhedral cone, i.e., $\set{C}$ is the conic hull of finitely many vectors. Combined with the constrained prior in \eqref{e:constrained_prior}, we can compute the normalization constant $K$ of the density of the prior conditioned on the face $\set{F}$ of $\set{C}$ as follows,
\begin{equation*}
    1 = \int_{\set{F}}\pi_{\vec{x}|\set{F}}(\vec{u}) \text{d}\vec{u}= K\int_{\set{F}}\exp\left(-\frac{\delta}{2}\|LF\vec{u}\|_2^2\right) \text{d}\vec{u} = K\delta^{-\text{dim}(\set{F})/2}\int_{\set{F}}\exp\left(-\frac{1}{2}\|LF\vec{v}\|_2^2\right) \text{d}\vec{v},
\end{equation*}
where we used that the face of any polyhedral cone is again a polyhedral cone, hence $c\set{F} = \set{F}$ for any $c > 0$ and we can take $\vec{x}_0 = \vec{0}$. The normalization constant $K$ is therefore
\begin{equation*}
    K = \frac{\delta^{\text{dim}(\set{F})/2}}{\int_{\set{F}}\exp\left(-\frac{1}{2}\|LF\vec{v}\|_2^2\right) \text{d}\vec{v}},
\end{equation*}
hence 
the distribution of the prior conditioned on the face satisfies
\begin{equation*}
    \pi_{\vec{x}| \set{F}}(\vec{u}) \,\propto\, \delta^{\text{dim}(\set{F})/2}\exp\left(-\frac{\delta}{2}\|LF\vec{u}\|_2^2\right),
\end{equation*}
where the proportionality does not depend on $\delta$ anymore.

Now we can obtain the distribution of the (hyper)parameters using Bayes' formula,
\begin{align*}
\pi_{\vec{x},\lambda,\delta | \vec{b},\set{F}}(\vec{u}, \lambda, \delta) 
\,\propto\,&\ \pi_{\vec{b}|\vec{x},\lambda,\delta}(F\vec{u})\pi_{\vec{x}| \delta, \set{F}}(\vec{u)} \pi_{\lambda}(\lambda) \pi_{\delta}(\delta)\\
\,\propto\,&\ \lambda^{m/2+\alpha_{\lambda} - 1}\delta^{\text{dim}(\set{F}(\vec{x}))/2+\alpha_{\delta} - 1}\\
       &\times\exp\left(-\frac{\lambda}{2}\|AF\vec{u} - \vec{b}\|_2^2 - \frac{\delta}{2}\|LF\vec{u}\|_2^2 - \beta_{\lambda}\lambda - \beta_{\delta}\delta\right),
\end{align*}
where $\set{F}(\vec{x})$ is the smallest face of $\set{C}$ that contains $\vec{x}$. Therefore, we can conclude the following conditional distributions of the hyperparameters,
\begin{align*}
    \nonumber
    \pi(\lambda\, |\, \vec{x}, \vec{b}, \set{F}) &\,\propto\, \lambda^{m/2 + \alpha_{\lambda}-1}\exp(-\beta_{\lambda}\lambda - \frac{\lambda}{2}\|A\vec{x} - \vec{b}\|_2^2), \text{ for } \lambda > 0\quad \text{and}\\
    \pi(\delta\, |\, \vec{x}, \vec{b}, \set{F}) &\,\propto\, \delta^{\text{dim}(\set{F}(\vec{x}))/2 + \alpha_{\delta}-1}\exp(-\beta_{\delta}\delta - \frac{\delta}{2}\|L\vec{x}\|_2^2), \text{ for } \delta > 0,
\end{align*}
or equivalently
\begin{align}
    \nonumber
    \lambda\, |\, \vec{x},\vec{b},\set{F} &\sim \Gamma(m/2 + \alpha_{\lambda}, \frac{1}{2}\|A\vec{x}-\vec{b}\|_2^2 + \beta_\lambda),\\
    \delta\, |\, \vec{x},\vec{b},\set{F} &\sim \Gamma(\text{dim}(\set{F}(\vec{x}))/2 + \alpha_{\delta}, \frac{1}{2}\|L\vec{x}\|_2^2 + \beta_\delta).
    \label{eq:delta_constrained}
\end{align}
From these conditional distributions, we obtain the Polyhedral Cone Hierarchical Gibbs Sampler \ref{sampler:PCHGS}.

\begin{algo}

\centering
\begin{minipage}{.6\textwidth}
\begin{algorithmic}
\STATE{ \textbf{Input:} $\vec{x}^0,\alpha_{\lambda}, \beta_{\lambda}, \alpha_{\delta}, \alpha_{\lambda}, k_{\max}$}

\FOR{$k = 1$ \TO $k_{\max}$ }

\STATE {Compute $(\lambda_k, \delta_k) \sim \pi_{\lambda, \delta|\vec{x},\vec{b}}$ as follows:} 
\STATE {\quad$\lambda_k \sim \Gamma\left(m/2+\alpha_{\lambda}, \frac12\|A\vec{x}^{k-1}-\vec{b}\|_2^2 + \beta_{\lambda}\right),$} 
\STATE {\quad$\delta_k \sim \Gamma\left(\text{dim}(\set{F}(\vec{x}^{k-1}))/2+\alpha_{\delta}, \frac12 \|L\vec{x}^{k-1}\|_2^2 + \beta_{\delta}\right).$} 
\STATE {Compute $\vec{x}^k \sim \pi_{\vec{x}|\vec{b},\lambda^k,\delta^k}$ using \eqref{e_randomized_constrained_least_squares_simple}} 

\ENDFOR

\RETURN{$\{(\vec{x}^k,\lambda_k,\delta_k)\}_{k=1,\dots, k_{\max}}$}

\end{algorithmic}
\end{minipage}
\caption{Polyhedral Cone Hierarchical Gibbs Sampler}
\label{sampler:PCHGS}
\end{algo}

The main difference with the ordinary Hierarchical Gibbs sampler is $\delta$, see the difference between \eqref{eq:delta_unconstrained} and \eqref{eq:delta_constrained}. It can be quite expensive to compute $\text{dim}(\set{F}(\vec{x}^k))$, but in a few cases there are simpler expressions. In the simplest setting where $\set{C} = \mathbb{R}^n$, i.e., the unconstrained setting, then the only face is the whole space $\mathbb{R}^n$, hence $\text{dim}(\set{F}(\vec{x}^k)) = n$ and the Polyhedral Cone Hierarchical Gibbs Sampler reduces to the ordinary Hierarchical Gibbs Sampler. In the more complicated setting where $\set{C} = \mathbb{R}^n_{\geq 0}$, i.e., nonnegativity constraints, then the faces are characterized by the zero values of the vector and $\text{dim}(\set{F}(\vec{x}))$ simplifies to the number of non-zero values elements of $\vec{x}$. The Polyhedral Cone Hierarchical Gibbs Sampler above then simplifies to the Nonnegative Hierarchical Gibbs Sampler of \cite{bardsley2020mcmc}.

\newpage
\section{Numerical examples}\label{sec_num_examples}
In this section, we present two numerical examples. First, we consider a one-dimensional deblurring problem and investigate the effect of constraints on the posterior distribution. Second, we consider a Bayesian hierarchical model for a two-dimensional computed tomography (CT) problem. For this problem, we use the Gibbs sampler as described in the Subsection \ref{ss_Gibbs} and take particular interest in efficiently and approximately solving the constrained linear least squares problems required for sampling.

\subsection{One-dimensional deblurring with different constraints}

Let us first consider a one-dimensional Gaussian deblurring problem defined by
\begin{equation}\label{eq:deblur_model}
    \vec{b} = A\vec{x} + \vec{e},
\end{equation}
for a true signal $\vec{x}\in [0,1]^n$, noise $\vec{e} \sim \mathcal{N}(\vec{0}, \lambda^{-1}I)$ with hyperparameter $\lambda > 0$ and forward operator $A$ defined by the Toeplitz matrix
\begin{equation*}
    A_{ij} = \frac{h}{\gamma\sqrt{2\pi}}\exp\left(-\frac12 \left(\frac{h(i-j)}{\gamma}\right)^2\right),\quad \text{ for } i,j = 1,\dots, n
\end{equation*}
where $n=128$, $h=1/n$ and $\gamma = 0.02$. Assume a priori that $\vec{x} \sim \mathcal{N}(\vec{0}, (\delta L^TL)^{-1})$, 
where $\delta > 0$ is a hyperparameter and $L$ is a periodic first-order finite difference matrix. Let the hyperparameters be fixed to $\lambda = 1000$ and $\delta = 150$.

\begin{figure}
    \centering
    \includegraphics[width = 0.8\textwidth]{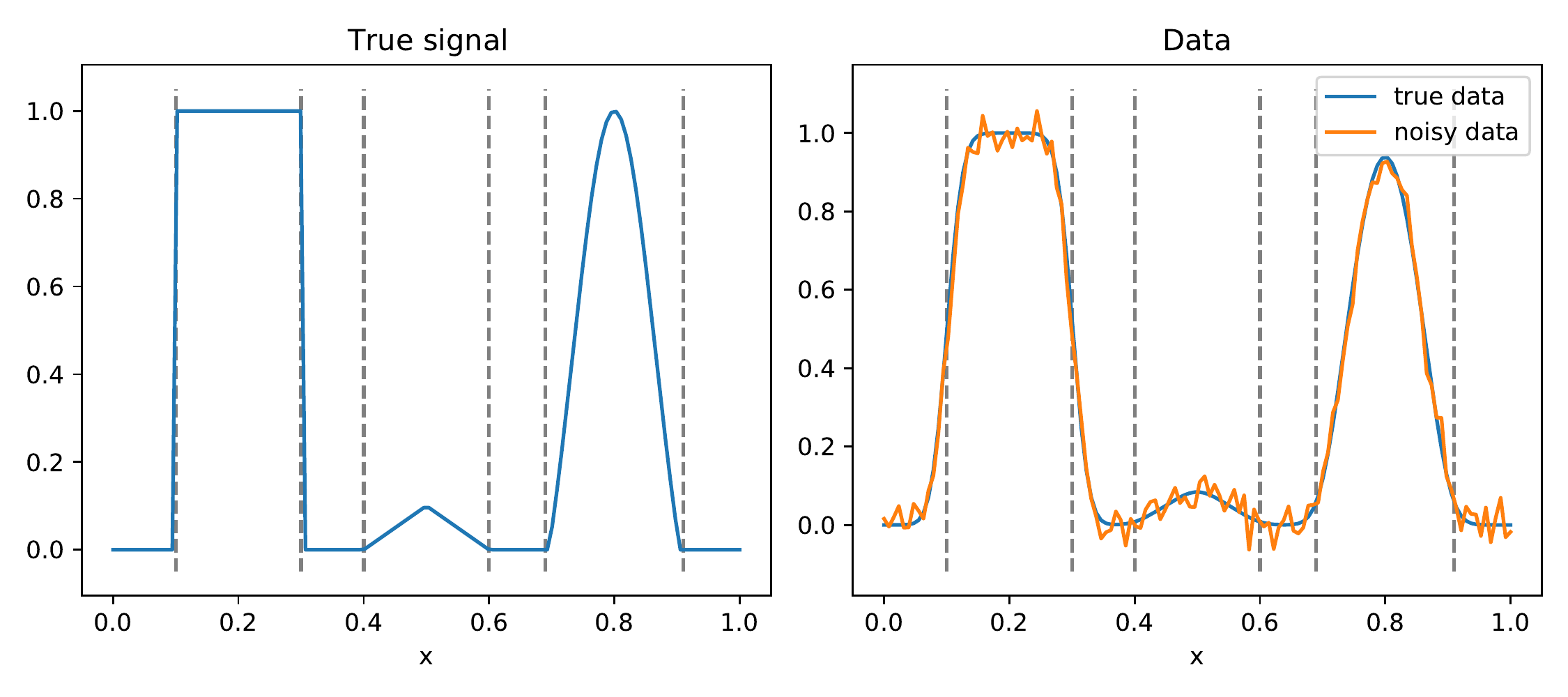}
    \caption{True signal, true data and noisy data for a Gaussian deblurring problem.}
    \label{fig:denoising_data}
\end{figure}
Figure \ref{fig:denoising_data} shows the true signal and the noisy measurements obtained through \eqref{eq:deblur_model} with $\lambda = 1000$. For the specific instance of $\vec{e}$ we have $\|\vec{e}\|/\|A\vec{x}\| \approx 6\%$. The true signal is divided into multiple components in order to illustrate the impact of constraints on different signal behaviour. These components include instantaneous changes between the extreme values $0$ and $1$, a small deviation from an extreme value and a smooth transition between extreme values. A lot of the components of the true signal have extreme values, hence we will be a priori interested in signals on the boundary of $[0,1]^n$, i.e., the set of signals for which at least one element is $0$ or $1$.

\begin{figure}
     \centering
     \begin{subfigure}[b]{\textwidth}
         \centering
         \includegraphics[width=\textwidth]{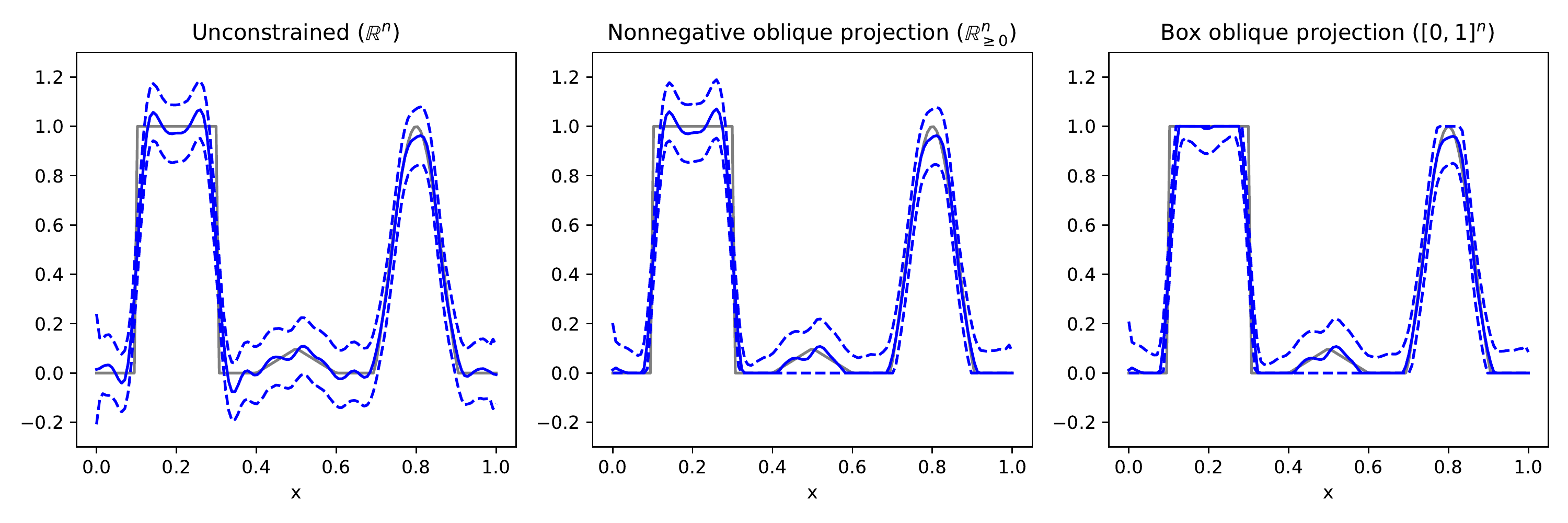}
         \caption{Different constraints with oblique projection.}
         \label{fig:denoising_comparison_oblique}
     \end{subfigure}
     \begin{subfigure}[b]{\textwidth}
         \raggedleft
         \includegraphics[width=0.67\textwidth]{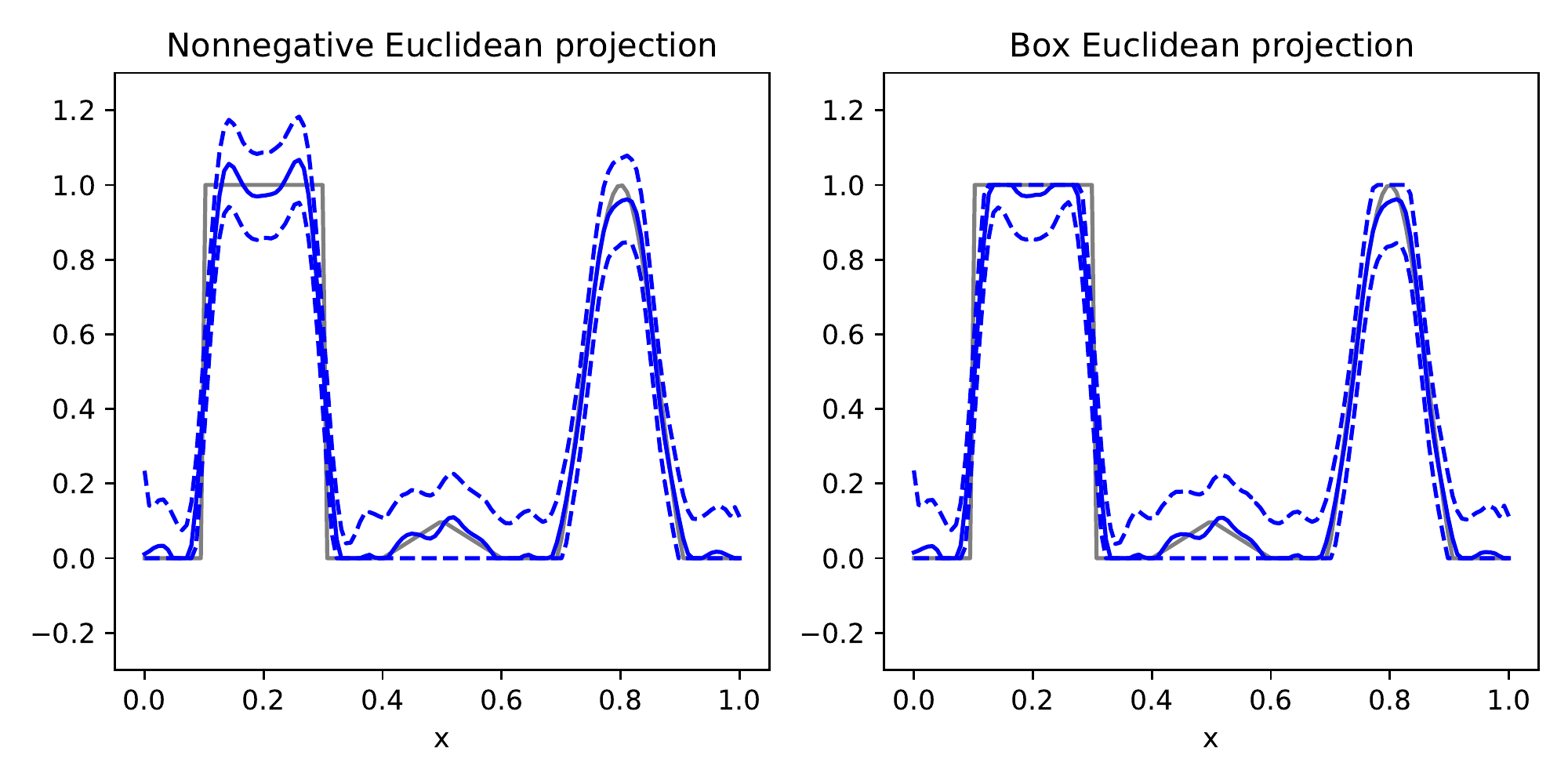}
         \caption{Different constraints with standard Euclidean projection.}
         \label{fig:denoising_comparison_projected}
     \end{subfigure}
        \caption{Component-wise median and $95\%$ credible intervals for 10000 samples of the Gaussian deblurring model.}
        \label{fig:denoising_comparison}
\end{figure}

Figure \ref{fig:denoising_comparison_oblique} shows the component-wise median and $95\%$ credible intervals for the posterior from $10000$ samples obtained by repeatedly solving optimization problem \eqref{e_randomized_constrained_least_squares_simple} in the unconstrained ($\reals^n$), nonnegative constrained ($\reals_{\geq 0}^n$) and box constrained ($[0,1]^n$) settings. Note that the fluctuating features near the extreme values $0$ and $1$ get flattened when applying the constraints. However, when applying the Euclidean projector to unconstrained samples, as illustrated in Figure \ref{fig:denoising_comparison_projected}, these fluctuations around the extreme values are stronger.

\begin{figure}
    \centering
    \includegraphics[width = 1.0\textwidth]{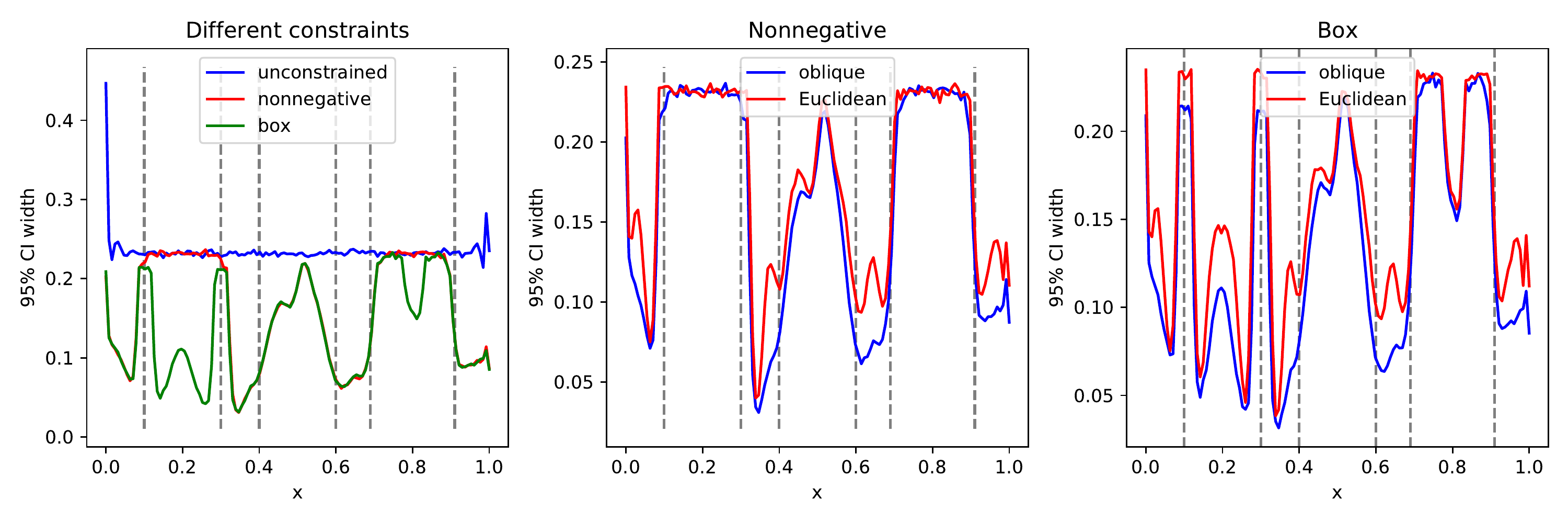}
    \caption{Component-wise width of $95\%$ credible intervals for different constraints and projectors (oblique and Euclidean) obtained from 10000 samples of the Gaussian deblurring model.}
    \label{fig:denoising_comparison_ci_width}
\end{figure}

As we use the median as central point estimate, we measured the variation of the samples using the width of component-wise credible intervals. Figure \ref{fig:denoising_comparison_ci_width} shows the width of $95\%$ component-wise credible intervals for the examples in Figure \ref{fig:denoising_comparison_oblique}. Figure \ref{fig:denoising_comparison_ci_width} shows that the width is reduced most where the signal values lie close the extreme values, but the width is also reduced for values close to extreme. Furthermore, using the oblique projection generally reduces the width more than the Euclidean projector. A possible explanation is that the Euclidean projection works component-wise, while the oblique projection takes interaction between components into account.

\subsection{Gibbs sampler for CT reconstruction with nonnegativity constraints}
\begin{figure}
    \centering
    \includegraphics[width = 0.9\textwidth]{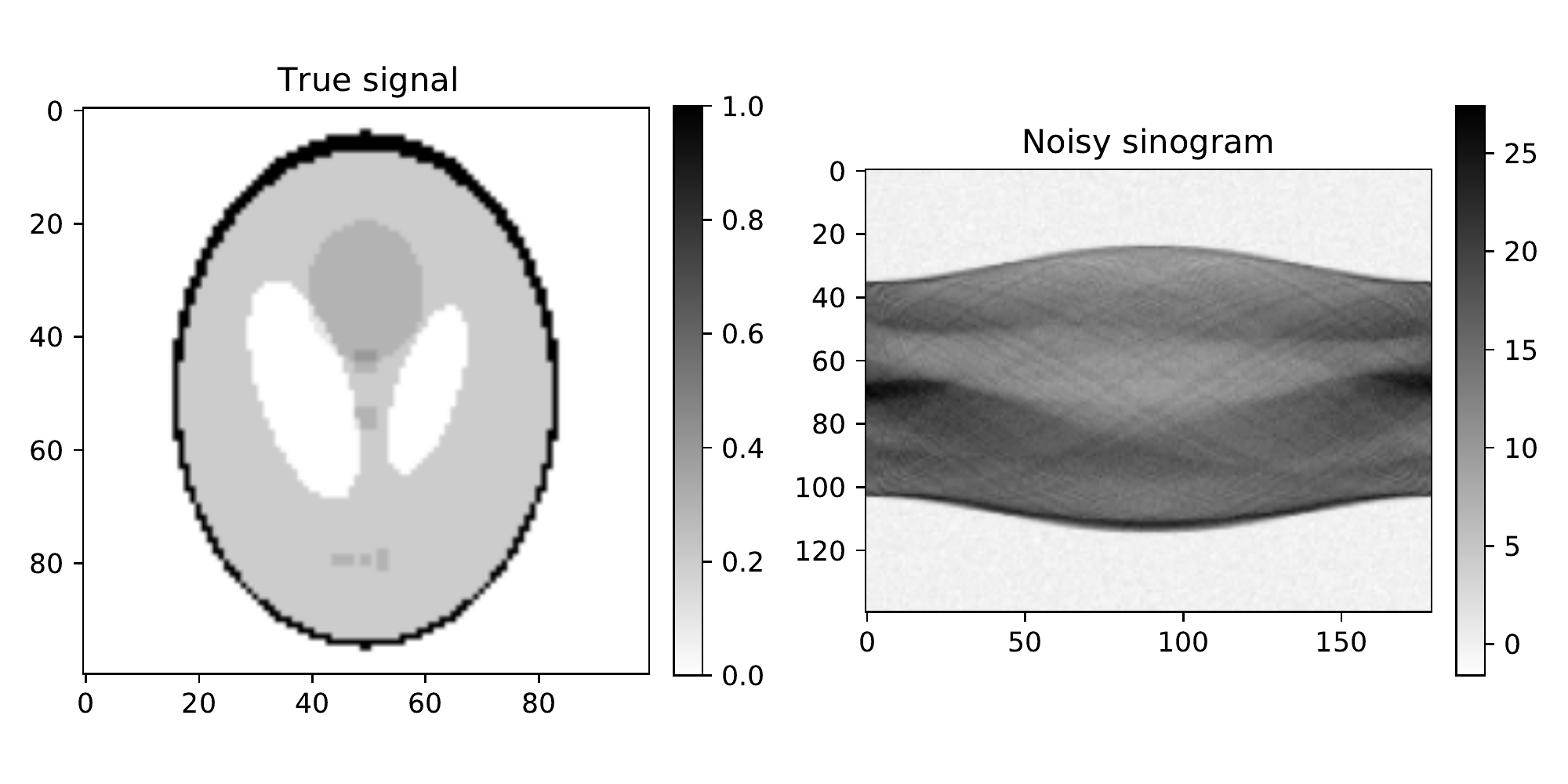}
    \caption{True Shepp-Logan phantom (left) and noisy sinogram (right).}
    \label{fig:CT_data}
\end{figure}

Let us now consider a CT problem \cite{hansen2021computed} of the form
\begin{equation*}
    \vec{b} = A\vec{x} + \vec{e},
\end{equation*}
where the true signal $\vec{x}\in \mathbb{R}^{100 \times 100}$ is the Shepp-Logan phantom shown in Figure \ref{fig:CT_data}, the noise satisfies $\vec{e} \sim \mathcal{N}(\vec{0}, \lambda^{-1}I)$ with hyperparameter $\lambda > 0$, and forward operator $A$ obtained from AIR Tools II \cite{hansen2018air} is a discretized Radon transform at $180$ angles with $140$ rays and a parallel-beam geometry. Assume a priori that $\vec{x} \sim \mathcal{N}(\vec{0}, (\delta L^TL)^{-1})$ with 
where $\delta > 0$ is a hyperparameter and $L$ is a  first-order finite difference matrix. We generated noisy data using $\lambda = 5$. For the specific instance of $\vec{e}$ we have $\|\vec{e}\|/\|A\vec{x}\| \approx 4\%$.

In CT, the signal $\vec{x}$ represents attenuation coefficients that are bounded from below by the corresponding background value, generally air. For simplicity, this constraint has been modelled as nonnegatity. In this experiment, we used the Polyhedral Cone Hierarchical Gibbs Sampler \label{sampler_PCHGS} to sample from the posterior both with and without nonnegativity constraints. The parameters for the hyperprior were chosen to be $\alpha_{\lambda} = \alpha_{\delta} = 1$ and $\beta_{\lambda} = \beta_{\delta} = 10^{-4}$, similar to \cite{bardsley2020mcmc}.

Repeatedly solving optimization problem \eqref{e_randomized_constrained_least_squares_simple} to high accuracy is a very costly procedure. Therefore, instead of solving the optimization problem from scratch at each iteration of the Gibbs sampler, we solve the optimization problem each time using a few iterations of FISTA \cite{beck2009fast} with the previous sample as warm-start. Using a small number of iterations, greatly reduces the computation time, but the samples become more correlated and do not have to be samples from the target distribution anymore. This is similar to using a few steps of CG in a gradient scan Gibbs sampler \cite[Algorithm 5.4]{bardsley2018computational},although this guarantees convergence to the target distribution.

We chose to run $100$ iterations of FISTA for each sample with the dynamic stepsize \linebreak $0.99(\lambda \|A^TA\|_2 + \delta \|L^TL\|_2)^{-1}$. This stepsize is a positive lower-bound on the inverse of the Lipschitz constant of the objective function and therefore guarantees convergence. We ran the algorithm for $15000$ samples and removed the first $1000$ samples as burn-in.

\begin{figure}
    \centering
    \includegraphics[width = 0.8\textwidth]{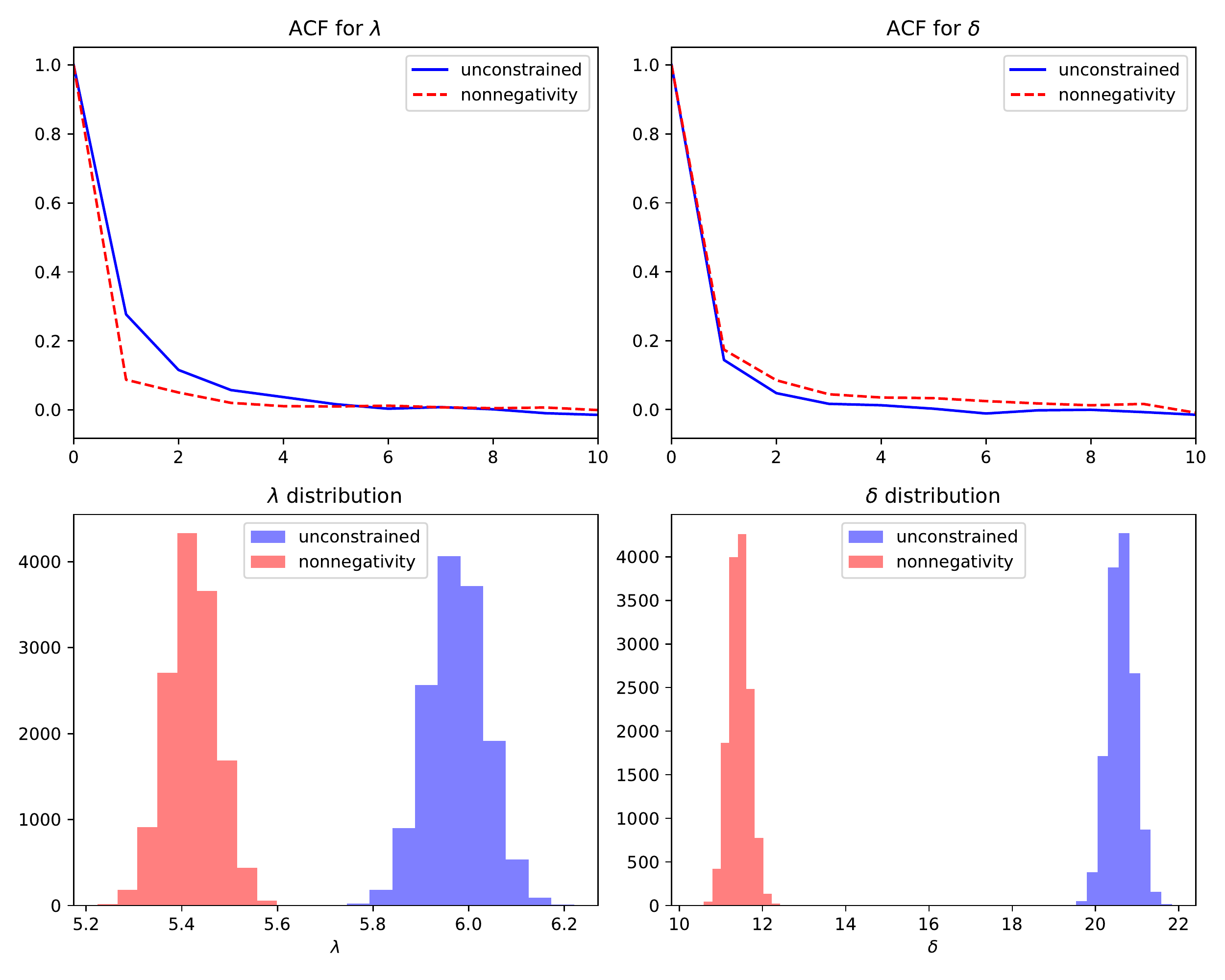}
    \caption{Autocorrelation functions (ACF) for the hyperparameters and hyperparameter distributions for unconstrained and nonnegative setting.}
    \label{fig:CT_hyperchains}
\end{figure}

Figure \ref{fig:CT_hyperchains} shows the autocorrelation function (ACF) and the distribution of the hyperparameters $\lambda$ and $\delta$. First, note that the autocorrelation function for $\lambda$ decays slightly faster in the nonnegative setting than the unconstrained setting. This shows that the nonnegative samples are slightly less correlated. Second, although the actual noise level is the same in both settings, the noise parameter $\lambda$ and the prior hyperparameter $\delta$ are noticeably smaller.

\begin{figure}
    \centering
    \includegraphics[width = 0.8\textwidth]{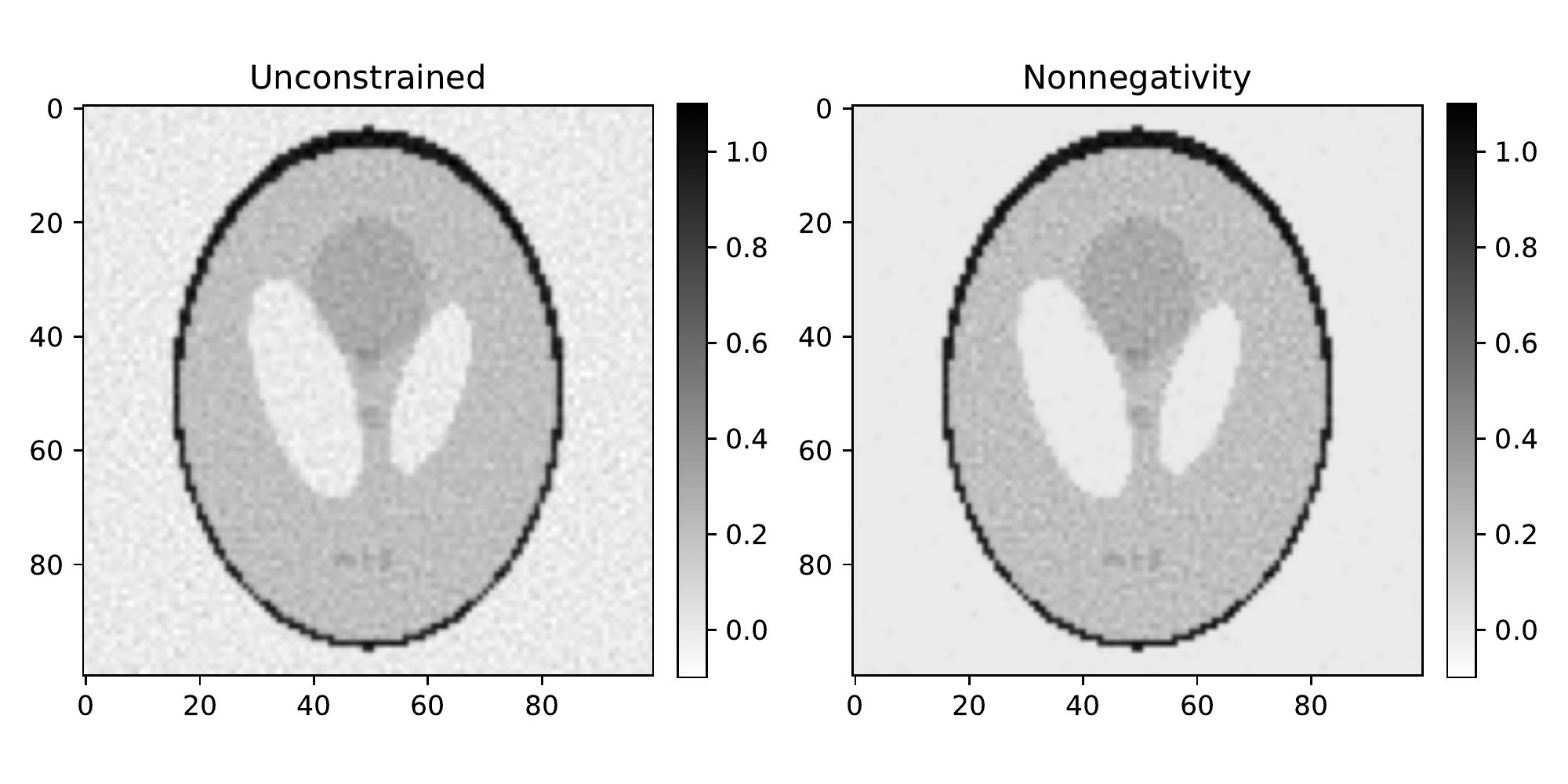}
    \caption{ Component-wise median for the unconstrained (left) and nonnegative (right) settings. Note that the range is the same.}
    \label{fig:CT_median}
\end{figure}

Figure \ref{fig:CT_median} shows the component-wise median for both the unconstrained and nonnegative setting. The main difference between the two medians lies in the background of the phantom. The unconstrained median has a non-uniform background containing a lot of small artefacts, while the nonnegative median is, besides a few pixels, uniformly zero in the background.

\begin{figure}
    \centering
    \includegraphics[width = 1.0\textwidth]{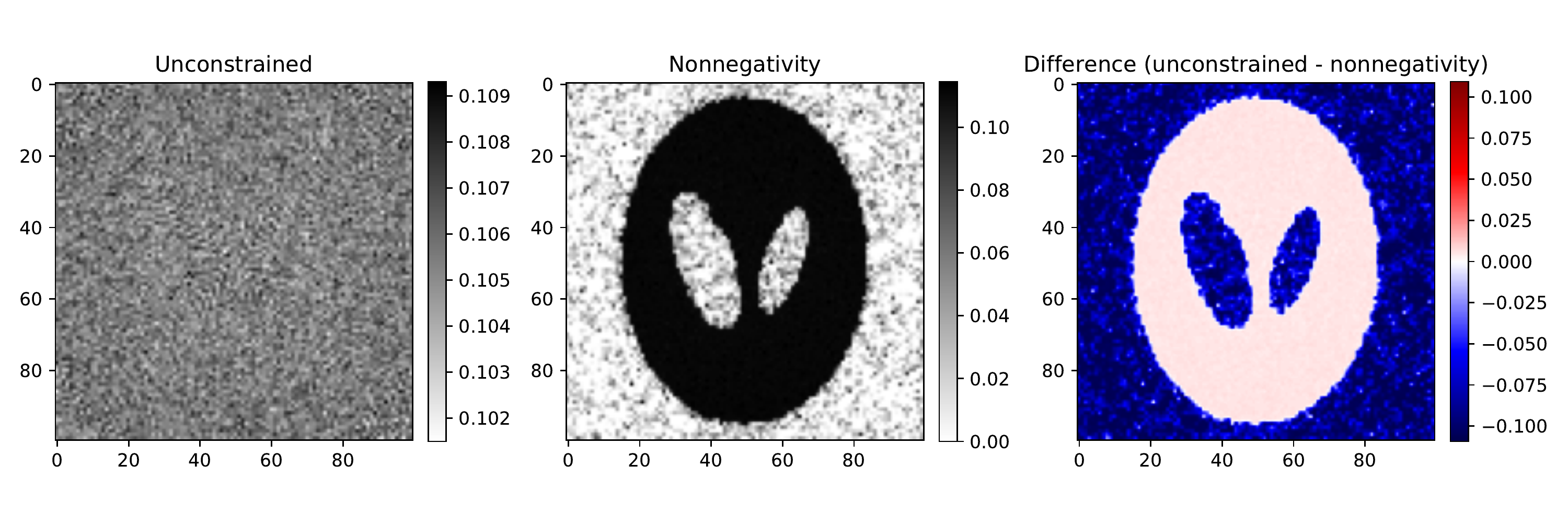}
    \caption{Component-wise width of $95\%$ credible intervals for the unconstrained (left) and nonnegative (middle) settings, together with their difference (right).}
    \label{fig:CT_CI_width}
\end{figure}

Another difference between the unconstrained and nonnegative setting is the component-wise variation. Figure \ref{fig:CT_CI_width} shows the width of the component-wise $95\%$ credible intervals together with the difference between the two settings. Note the different range of the unconstrained and nonnegative settings. Just like in the deblurring experiment, the width is greatly reduced in the components that are close to the extreme values. Therefore, the uncertainty can be used to easily distinguish between the actual object and the background. The uncertainty inside the object is slightly larger with nonnegativity constraints, which is due to the smaller $\lambda$ and $\delta$ hyperparameters as observed in Figure \ref{fig:CT_hyperchains}.

\section{Conclusion}

We have presented a method for handling constraints in Bayesian inference that puts positive probability on the boundary of the constraint set. In general, the method works by projecting posterior samples from outside the constraint set onto the constraint set. Therefore, the method can be interpreted as post-processing the posterior by projecting the density onto the constraint set.

The posterior sampling and projection steps can be combined if the posterior distribution is Gaussian and the projection is with respect to the posterior precision matrix. Samples from such a post-processed posterior can be obtained by solving perturbed constrained quadratic optimization problems. Although this distribution is difficult to describe in closed form for general constraint sets, we were able to characterize the distribution when the constraint set is a polyhedral set. We have proven that the projected posterior on a polyhedral set consists of densities on each of the faces and these densities are proportional to the unprojected posterior density.

To apply the theory, we considered Bayesian linear inverse problems. For such problems, sampling from the projected posterior can be achieved by solving perturbed constrained linear least squares problems. Furthermore, we considered a Bayesian hierarchical model and derived a Gibbs sampler for the constrained problem when the constraint set is a polyhedral cone.

We tested the projection method on deblurring and CT test cases, for which component-wise bounds are natural constraints. These numerical experiments have shown that the projected posterior can greatly reduce the uncertainty of vector components that are close to their bounds. Furthermore, the experiments suggest that using the oblique projection instead of the Euclidean projection gives better results.

One major issue with sampling by means of solving perturbed constrained least squares problems is the computational cost. The computational cost of solving these optimization problems accurately can be very high, yet it is still unknown to what extent accurate solutions are necessary. Therefore, further research should focus on determining what the effect that inaccurately solving the optimization problems has on the distribution of the samples.

Another topic of further study is identifying more techniques for efficient post-processing of posteriors beyond constraints. For example, adding penalization functions to the randomized least squares problems results in new modified posteriors. The computational cost of sampling from this new posterior is similar to the projected posterior, but it can introduce different regularization-like effects. 
\newpage
\bibliographystyle{siamplain}
\bibliography{references}

\appendix
\section{Miscellaneous proofs and computations}
\subsection{Linear algebra lemma}\label{subsec:appendix_proof}
\begin{lemma}\label{lem:basis_transform}
If $\vec{u}_1, \dots, \vec{u}_k, \vec{u}_{k+1}, \dots, \vec{u}_n \in \mathbb{R}^n$ is an orthonormal basis and $\Sigma$ is a positive definite matrix, then $\vec{u}_1, \dots, \vec{u}_k, \Sigma \vec{u}_{k+1}, \dots, \Sigma \vec{u}_n \in \mathbb{R}^n$ is a basis for $\mathbb{R}^n$.
\end{lemma}
\begin{proof}
Let $U = [U_1, U_2]$ where $U_1$ and $U_2$ are matrices with the first $k$ and last $n-k$ vectors $u_i$ as columns respectively. We will proof the claim by showing that $\hat{U} = [U_1, \Sigma U_2]$ has full rank. Note that
$$
U^T\hat{U} = 
\begin{bmatrix}
I & U_1^T\Sigma U_2\\
0 & U_2^T\Sigma U_2
\end{bmatrix},
$$
has full rank, because $U_2^T\Sigma U_2$ is positive definite. Therefore, $\hat{U} = U(U^T\hat{U})$ also has full rank.

\end{proof}

\subsection{Analytic computations}\label{subsec:appendix_analytical_examples}
\subsubsection{Halfspace}
Define the halfspace $\set{C} = \{\vec{x}\in \reals^m \,|\, \vec{a}^T\vec{x} \leq b\}$, where $\vec{a} \in \reals^n$ is a nonzero normal vector and $b \in \reals$. Denote by $F \in \reals^{n\times(n-1)}$ a matrix whose columns form an orthonormal basis for $\nullspace(\vec{a}^T)$, then the boundary of the halfspace defined can be parameterized by $\vec{x}_0 + F\vec{u}$ for $\vec{u}\in \reals^{n-1}$. Let $\vec{x}^\star \sim \mathcal{N}(\bm{\mu}, \Sigma)$ and $\set{E} \subseteq \bd(\set{C})$ be measurable, then
\begin{align*}
    \mathbb{P}\left(\Pi_{\set{C}}^{\Sigma^{-1}}(\vec{x}^\star) \in \set{E}\right) &= |\det\left(\begin{bmatrix}\Sigma \vec{a} & F\end{bmatrix}\right)|
    \int_{\reals^{n-1}}\int_{0}^\infty\pi_{\vec{x}^{\star}}(\vec{x}_0+F\vec{u}+t\Sigma \vec{a}) \text{d}t \, \text{d}\vec{u}\\
    &=\int_{\reals^{n-1}}\frac{\vec{a}^T\Sigma \vec{a}}{\|\vec{a}\|_2}\int_{0}^\infty\pi_{\vec{x}^{\star}}(\vec{x}_0+F\vec{u}+t\Sigma \vec{a}) \text{d}t \, \text{d}\vec{u},
\end{align*}
hence the density on the boundary can be written as
\begin{equation*}
    \pi_{\bd}(\vec{u}) = \frac{\vec{a}^T\Sigma \vec{a}}{\|\vec{a}\|_2}\int_{0}^\infty\pi_{\vec{x}^{\star}}(\vec{x}_0+F\vec{u}+t\Sigma \vec{a}) \text{d}t.
\end{equation*}

The integral can be computed using the following integral identity \cite{jeffrey2008handbook}: for constants $a < 0$ and $b,c \in \mathbb{R}$,
\[
\int_0^\infty \exp(at^2 + bt + c) \text{d}t = \frac{\sqrt{\pi}}{2\sqrt{-a}}\exp\left(-\frac{b^2}{4a}+c\right)\erfc\left(-\frac{b}{2\sqrt{-a}}\right).
\]

\subsubsection{Disc}
Suppose that $\set{C}$ is a unit disc defined by $\set{C} = \{\vec{x} \in \reals^2 \,|\, \|\vec{x}\|_2 \leq 1\}$ with boundary parameterization $\vec{n}(u) := (\cos(u), \sin(u))^T$ for $u \in [0, 2\pi)$. For the two-dimensional ball, we can derive an exact distribution on its boundary in a similar way to the boundary of a halfspace. Let $\set{E} \subseteq \bd(\set{C})$ be measurable, then
\[
 \mathbb{P}\left(\Pi_{\set{C}}^{\Sigma^{-1}}(\vec{x}^\star) \in \set{E}\right) = \int_\set{E} \int_0^\infty \pi_{\vec{x}^\star}(\vec{n}(u) + t\Sigma \vec{n}(u)) |J(t, u)| \text{d}t\,\text{d}u,
\]
where 
\[
|J(t, \theta)| = \det \begin{bmatrix}
  \Sigma \vec{n}(u) & (I+t\Sigma)R\vec{n}(u)
\end{bmatrix}, \quad \text{with}\quad R=\begin{bmatrix}
  0 & -1 \\ 1 & 0
\end{bmatrix},
\]
or equivalently
\[
|J(t, \theta)| = \det \begin{bmatrix}
  \Sigma \vec{n}(u) & R\vec{n}(u)
\end{bmatrix} + t \det (\Sigma) =: K(u) + t \det (\Sigma).
\]
Therefore, the resulting boundary distribution is given by
\[
\pi_{\bd}(u) = \int_0^\infty (K(u) + t \det (\Sigma))\pi_{\vec{x}^\star}(\vec{n}(u) + t\Sigma \vec{n}(u)) \text{d}t.
\]

The integral can be computed using the following integral identity \cite{jeffrey2008handbook}: for constants $a > 0$ and $b,c,d,f \in \mathbb{R}$,
\[
\int_0^\infty (d+ft)\exp(-\frac{1}{2}(at^2 + bt + c)) \text{d}t = \frac{e^{-c/2}}{4a^{3/2}}\left(4f\sqrt{a} + \sqrt{2\pi}(2ad-bf)\exp\left(\frac{b^2}{8a}\right)\erfc\left(\frac{b}{\sqrt{8a}}\right)\right).
\]

\end{document}